\documentclass[twoside]{amsart}
\usepackage{amsfonts,amsthm}
\usepackage{enumerate}
\usepackage{amsmath}
\usepackage{amssymb}
\usepackage{amsxtra}
\usepackage{latexsym}

\theoremstyle{plain}
\newtheorem*{thm A}{Theorem~A}
\newtheorem*{thm B}{Theorem~B}
\newtheorem*{thm C}{Theorem~C}

\newtheorem*{main 1}{Theorem~1}
\newtheorem*{main 2}{Theorem~2}
\newtheorem*{main 3}{Theorem~3}
\newtheorem*{coro 1}{Corollary~1}
\newtheorem*{coro 2}{Corollary~2}
\newtheorem*{lem 1}{Lemma~1}
\newtheorem*{lem 2}{Lemma~2}
\newtheorem*{lem 3}{Lemma~3}
\newtheorem*{pro A}{Proposition~A}
\newtheorem*{pro B}{Proposition~B}
\newtheorem*{rem}{Remark}

\newtheorem{theorem}{Theorem}[section]
\newtheorem{corollary}[theorem]{Corollary}
\newtheorem{lemma}[theorem]{Lemma}
\newtheorem{proposition}[theorem]{Proposition}

\begin{document}

\title[Isometric Reeb flow]{Isometric Reeb flow in complex hyperbolic quadrics}
\vspace{0.2in}

\author[Young Jin Suh]{Young Jin Suh}
\address{\newline  Young Jin Suh
\newline Department of Mathematics,
\newline College of Natural Sciences
\newline Kyungpook National University,
\newline Daegu 702-701, Republic of Korea}
\email{yjsuh@knu.ac.kr}

\footnotetext[1]{{\it 2010 Mathematics Subject Classification} : Primary 53C40; Secondary 53C55, 53D15.}
\footnotetext[2]{{\it Key words} : Real hypersurface; Reeb flow; complex hyperbolic quadric.}
\thanks{* This work was supported by Grant Proj. No. NRF-2015-R1A2A1A-01002459 from National Research Foundation, Korea.}

\begin{abstract}
We classify real hypersurfaces with isometric Reeb flow in the complex hyperbolic quadrics ${Q^*}^{m} = SO^{o}_{2,m}/SO_mSO_2$, $m \geq 3$. We show that $m$ is even, say $m = 2k$, and any such hypersurface becomes an open part of a  tube around a $k$-dimensional complex hyperbolic space ${\mathbb C}H^k$ which is embedded canonically in ${Q^*}^{2k}$ as a totally geodesic complex submanifold or a horosphere whose center at infinity is $\frak A$-isotropic singular. As a consequence of the result, we get the non-existence of real hypersurfaces with isometric Reeb flow in odd-dimensional complex quadrics ${Q^*}^{2k+1}$, $k \geq 1$.
\end{abstract}

\maketitle

\section{Introduction}
\setcounter{equation}{0}
\renewcommand{\theequation}{0.\arabic{equation}}
\vspace{0.13in}

Let $M$ be a real hypersurface in a K\"{a}hler manifold $\bar{M}$. The complex structure $J$ on $\bar{M}$ induces locally an almost contact metric structure $(\phi,\xi,\eta,g)$ on $M$. In the context of contact geometry, the unit vector field $\xi$ is often referred to as the Reeb vector field on $M$ and its flow is known as the Reeb flow. The Reeb flow has been of significant interest in recent years, for example in relation to the Weinstein Conjecture. We are interested in the Reeb flow in the context of Riemannian geometry, namely in the classification of real hypersurfaces with isometric Reeb flow in homogeneous K\"{a}hler manifolds.

For the complex hyperbolic space ${\mathbb C}H^m$ a full classification was obtained by Montiel and Romero in \cite{MR}. He proved that the Reeb flow on a real hypersurface in ${\mathbb C}H^m = SU_{1,m}/S(U_mU_1)$ is isometric if and only if $M$ is an open part of a tube around a totally geodesic ${\mathbb C}H^k \subset {\mathbb C}H^m$ for some $k \in \{0,\ldots,m-1\}$. For the complex hyperbolic $2$-plane Grassmannian $G_2^{*}({\mathbb C}^{m+2})= SU_{2,m}/S(U_mU_2)$ the classification was obtained by Suh in \cite{S5} and \cite{S6}. The Reeb flow on a real hypersurface in $G^{*}_2({\mathbb C}^{m+2})$ is isometric if and only if $M$ is an open part of a tube around a totally geodesic $G_2^{*}({\mathbb C}^{m+1}) \subset G_2^{*}({\mathbb C}^{m+2})$ or a horosphere with singular normal $JN{\in}{\mathfrak J}N$. In this paper we investigate this problem for the complex hyperbolic quadric ${Q^*}^m = SO_{2,m}/SO_mSO_2$. In view of the previous two results a natural expectation is that the classification involves at least the totally geodesic ${Q^*}^{m-1} \subset {Q^*}^m$. Surprisingly, this is not the case. Our main result states:

\begin{theorem} \label{mainresult}
Let $M$ be a real hypersurface of the complex hyperbolic quadric ${Q^*}^{m} = SO^{o}_{2,m}/SO_mSO_2$, $m \geq 3$. The Reeb flow on $M$ is isometric if and only if $m$ is even, say $m = 2k$, and $M$ is an open part of a tube around a totally geodesic ${\mathbb C}H^k \subset {Q^*}^{2k}$ or a horosphere whose center at infinity is $\mathfrak A$-isotropic singular.
\end{theorem}

Every tube around a totally geodesic ${\mathbb C}H^k \subset {Q^*}^{2k}$ is a homogeneous hypersurface. In fact, the closed subgroup $U_{1,k}$ of $SO^{o}_{2,2k}$ acts on ${Q^*}^{2k}$ with cohomogeneity one. The singular orbit is  totally geodesic ${\mathbb C}H^k \subset {Q^*}^{2k}$ and the principal orbits is the tubes around of this singular orbit. So as a corollary we get:

\begin{corollary}
Let $M$ be a connected complete real hypersurface in the complex quadric ${Q^*}^{2k}$, $k \geq 2$. If the Reeb flow on $M$ is isometric, then $M$ is a homogeneous hypersurface of ${Q^*}^{2k}$.
\end{corollary}

Our paper is organized as follows. In Section 2 we present basic material about the complex quadric ${Q^*}^m$, including its Riemannian curvature tensor and a description of its singular tangent vectors. Apart from the complex structure $J$
there is another distinguished geometric structure on ${Q^*}^m$, namely a parallel rank two vector bundle ${\mathfrak A}$ which contains an $S^1$-bundle of real structures on the tangent spaces of ${Q^*}^m$. This geometric structure determines a maximal ${\mathfrak A}$-invariant subbundle ${\mathcal Q}$ of the tangent bundle $TM$ of a real hypersurface $M$ in ${Q^*}^m$. In Section 3 we investigate the geometry of this subbundle ${\mathcal Q}$. In Section 4 we describe the canonical embedding of ${\mathbb C}H^k$ into ${Q^*}^{2k}$ as a totally geodesic complex submanifold and investigate the geometry of the tubes around ${\mathbb C}H^k$. We will show that the Reeb flow on these tubes is isometric. In Section 5 we determine some geometric consequences of the Codazzi equation for real hypersurfaces in ${Q^*}^m$.

Finally, in Section 6, we present the proof of Theorem \ref{mainresult}. The first step is to prove that the normal bundle of a real hypersurface $M$ with isometric Reeb flow consists of a particular type of singular tangent vectors of ${Q^*}^m$, the so-called ${\mathfrak A}$-isotropic vectors. The next step is to show that ${\mathcal Q}$ is invariant under the shape operator of $M$. Putting all this information into the Codazzi equation then allows us to compute explicitly the principal curvatures and principal curvatures spaces of $M$. A particular consequence of this is that $m$ is even, say $m = 2k$.
Using Jacobi field theory we then show that $M$ has a smooth focal manifold at a constant distance which is embedded in ${Q^*}^{2k}$ as a totally geodesic complex submanifold of complex dimension $k$.  Corresponding  to Klein's \cite{K} classification, if we use the classification theory of totally geodesic submanifolds in noncompact complex quadrics ${Q^*}^{2k}$ we will show that this focal manifold is a totally geodesic ${\mathbb C}H^k \subset {Q^*}^{2k}$.

\section{The complex hyperbolic quadric} \label{quadric}

Let us denote by $C_1^{m+2}$ an indefinite complex Euclidean space ${\mathbb C}^{m+2}$, on which the indefinite Hermitian product
$$H(z,w)=-z_1{\bar w}_1+z_2{\bar w}_2+{\cdots}+z_{n+2}{\bar w}_{n+2}$$
is negative definite.

The homogeneous quadratic equation $z_1^2 + \ldots + z_{m}^2 - z_{m+1}^2-z_{m+2}^2 = 0$
consists of the points in ${\mathbb C}_1^{m+2}$ defines a noncompact complex hyperbolic quadric ${Q^*}^m=SO^{o}_{2,m}/SO_2SO_m$ which can be immersed in the $(m+1)$-dimensional in complex hyperbolic space ${\mathbb C}H^{m+1} = SU_{1,m+1}/S(U_{m+1}U_1)$. The complex hypersurface ${Q^*}^m$ in ${\mathbb C}H^{m+1}$ is known as the $m$-dimensional complex hyperbolic quadric. The complex structure $J$ on ${\mathbb C}H^{m+1}$ naturally induces a complex structure on  ${Q^*}^m$ which we will denote by $J$ as well. We equip ${Q^*}^m$ with the Riemannian metric $g$ which is induced from the Begerman metric on ${\mathbb C}H^{m+1}$ with constant holomorphic sectional curvature $4$.  For $m \geq 2$ the triple $({Q^*}^m,J,g)$ is a Hermitian symmetric space of rank two and its minimal sectional curvature is equal to $-4$. The $1$-dimensional quadric ${Q^*}^1$ is isometric to the $2$-dimensional real hyperbolic space ${\mathbb R}H^2=SO^{o}_{1,2}/SO_1SO_2$.  The $2$-dimensional complex quadric ${Q^*}^2$ is isometric to the Riemannian product of complex hyperbolic spaces ${\mathbb C}H^1 \times {\mathbb C}H^1$. We will assume $m \geq 3$ for the main part of this paper.

For a nonzero vector $z \in {\mathbb C}_1^{m+2}$ we denote by $[z]$ the complex span of $z$, that is, $[z] = \{\lambda z \mid \lambda \in {\mathbb C}\} $.
Note that by definition $[z]$ is a point in ${\mathbb C}H^{m+1}$.
As usual, for each $[z] \in {\mathbb C}H^{m+1}$ we identify $T_{[z]}{\mathbb C}H^{m+1}$ with the orthogonal complement ${\mathbb C}_1^{m+2} \ominus [z]$ of $[z]$ in ${\mathbb C}_1^{m+2}$. For $[z] \in {Q^*}^m$ the tangent space $T_{[z]}{Q^*}^m$ can then be identified canonically with the orthogonal complement ${\mathbb C}_1^{m+2} \ominus ([z] \oplus [\bar{z}])$ of $[z] \oplus [\bar{z}]$ in ${\mathbb C}_1^{m+2}$. Note that $\bar{z} \in \nu_{[z]}{Q^*}^m$ is a unit normal vector of ${Q^*}^m$ in ${\mathbb C}H^{m+1}$ at the point $[z]$.

We denote by $A_{\bar{z}}$ the shape operator of ${Q^*}^m$ in ${\mathbb C}H^{m+1}$ with respect to $\bar{z}$. Then we have $A_{\bar{z}}w = \overline{w}$ for all $w \in T_{[z]}Q*^m$, that is,  $A_{\bar{z}}$ is just complex conjugation restricted to $T_{[z]}{Q^*}^m$.
The shape operator $A_{\bar{z}}$ is an antilinear involution on the complex vector space $T_{[z]}{Q^*}^m$ and
$$T_{[z]}{Q^*}^m = V(A_{\bar{z}}) \oplus JV(A_{\bar{z}}),$$
where $V(A_{\bar{z}}) =  {\mathbb R}_1^{m+2} \cap T_{[z]}{Q^*}^m$ is the $(+1)$-eigenspace and $JV(A_{\bar{z}}) = i{\mathbb R}_1^{m+2} \cap T_{[z]}{Q^*}^m$ is the $(-1)$-eigenspace of $A_{\bar{z}}$.  Geometrically this means that the shape operator $A_{\bar{z}}$ defines a real structure on the complex vector space $T_{[z]}{Q^*}^m$. Recall that a real structure on a complex vector space $V$ is by definition an antilinear involution $A : V \to V$. Since the normal space $\nu_{[z]} {Q^*}^m$ of ${Q^*}^m$ in ${\mathbb C}H_1^{m+1}$ at ${[z]}$ is a complex subspace of $T_{[z]}{\mathbb C}H^{m+1}$ of complex dimension one, every normal vector in $\nu_{[z]}{Q^*}^m$ can be written as $\lambda \bar{z}$ with some $\lambda \in  {\mathbb C}$. The shape operators $A_{\lambda\bar{z}}$ of ${Q^*}^m$ define a rank two vector subbundle ${\mathfrak A}$ of the endomorphism bundle ${\rm End}(T{Q^*}^m)$. Since the second fundamental form of the embedding ${Q^*}^m \subset {\mathbb C}H^{m+1}$ is parallel (see e.g.\ \cite{S}), ${\mathfrak A}$ is a parallel subbundle of ${\rm End}(T{Q^*}^m)$. For $\lambda \in S^1 \subset {\mathbb C}$ we again get a real structure $A_{\lambda\bar{z}}$ on $T_{[z]}{Q^*}^m$ and we have $V(A_{\lambda\bar{z}}) = \lambda V(A_{\bar{z}})$. We thus have an $S^1$-subbundle of ${\mathfrak A}$ consisting of real structures on the tangent spaces of ${Q^*}^m$.

The Gauss equation for the complex hypersurface ${Q^*}^m \subset {\mathbb C}H^{m+1}$ implies that the Riemannian curvature tensor $R$ of ${Q^*}^m$ can be expressed in terms of the Riemannian metric $g$, the complex structure $J$ and a generic real structure $A$ in ${\mathfrak A}$:
\begin{eqnarray*}
R(X,Y)Z & = &-g(Y,Z)X + g(X,Z)Y \\
& & - \, g(JY,Z)JX + g(JX,Z)JY + 2g(JX,Y)JZ \\
 & & -\, g(AY,Z)AX + g(AX,Z)AY \\
& & - \,g(JAY,Z)JAX + g(JAX,Z)JAY.
\end{eqnarray*}
Note that the complex structure $J$ anti-commutes with each endomorphism $A \in {\mathfrak A}$, that is, $AJ = -JA$.

A nonzero tangent vector $W \in T_{[z]}{Q^*}^m$ is called singular if it is tangent to more than one maximal flat in ${Q^*}^m$. There are two types of singular tangent vectors for the complex quadric ${Q^*}^m$:
\begin{enumerate}[\rm (i)]
\item If there exists a real structure $A \in {\mathfrak A}_{[z]}$ such that $W \in V(A)$, then $W$ is singular. Such a singular tangent vector is called ${\mathfrak A}$-principal.
\item If there exist a real structure $A \in {\mathfrak A}_{[z]}$ and orthonormal vectors $X,Y \in V(A)$ such that $W/||W|| = (X+JY)/\sqrt{2}$, then $W$ is singular. Such a singular tangent vector is called ${\mathfrak A}$-isotropic.
\end{enumerate}
Basic complex linear algebra shows that for every unit tangent vector $W \in T_{[z]}{Q^*}^m$ there exist a real structure $A \in {\mathfrak A}_{[z]}$ and orthonormal vectors $X,Y \in V(A)$ such that
\[
W = \cos(t)X + \sin(t)JY
\]
for some $t \in [0,\pi/4]$. The singular tangent vectors correspond to the values $t = 0$ and $t = \pi/4$.

\section{The maximal ${\mathfrak A}$-invariant subbundle ${\mathcal Q}$ of $TM$}

Let $M$ be a  real hypersurface in ${Q^*}^m$ and denote by $(\phi,\xi,\eta,g)$ the induced almost contact metric structure on $M$ and by $\nabla$ the induced Riemannian connection on $M$. Note that $\xi = -JN$, where $N$ is a (local) unit normal vector field of $M$. The vector field $\xi$ is known as the Reeb vector field of $M$. If the integral curves of $\xi$ are geodesics in $M$, the hypersurface $M$ is called a Hopf hypersurface. The integral curves of $\xi$ are geodesics in $M$ if and only if $\xi$ is a principal curvature vector of $M$ everywhere.
The tangent bundle $TM$ of $M$ splits orthogonally into  $TM = {\mathcal C} \oplus {\mathcal F}$, where ${\mathcal C} = {\rm ker}(\eta)$ is the maximal complex subbundle of $TM$ and ${\mathcal F} = {\mathbb R}\xi$. The structure tensor field $\phi$ restricted to ${\mathcal C}$ coincides with the complex structure $J$ restricted to ${\mathcal C}$, and we have $\phi \xi = 0$. We denote by $\nu M$ the normal bundle of $M$.

We first introduce some notations. For a fixed real structure $A \in {\mathfrak A}_{[z]}$ and $X \in T_{[z]}M$ we decompose $AX$ into its tangential and normal component, that is,
\[
AX = BX + \rho(X)N
\]
where $BX$ is the tangential component of $AX$ and
\[
\rho(X) = g(AX,N) = g(X,AN) = g(X,AJ\xi) = g(JX,A\xi).
\]
Since $JX = \phi X + \eta(X)N$ and $A\xi = B\xi + \rho(\xi)N$ we also have
\[
\rho(X) = g(\phi X,B\xi) + \eta(X)\rho(\xi)  = \eta(B\phi X) + \eta(X)\rho(\xi).
\]
We also define
\[
\delta = g(N,AN) = g(JN,JAN) = -g(JN,AJN) = -g(\xi,A\xi).
\]

At each point $[z] \in M$ we define
\[
{\mathcal Q}_{[z]} = \{X \in T_{[z]}M \mid AX \in T_{[z]}M\ {\rm for\ all}\ A \in {\mathfrak A}_{[z]}\},
\]
which is the maximal ${\mathfrak A}_{[z]}$-invariant subspace of $T_{[z]}M$. Then by using the same method for real hypersurfaces in ${Q^*}^m$ as in Suh \cite{S7} and \cite{S9} we get the following

\begin{lemma} \label{principal}
Let $M$ be a real hypersurface in ${Q^*}^m$. Then the following statements are equivalent:
\begin{enumerate}[\rm (i)]
\item The normal vector $N_{[z]}$ of $M$ is ${\mathfrak A}$-principal,
\item  ${\mathcal Q}_{[z]} = {\mathcal C}_{[z]}$,
\item  There exists a real structure $A \in {\mathfrak A}_{[z]}$ such that $AN_{[z]} \in {\mathbb C}\nu_{[z]}M$.
\end{enumerate}
\end{lemma}

\medskip
Assume now that the normal vector $N_{[z]}$ of $M$ is not ${\mathfrak A}$-principal.
Then there exists a real structure $A \in {\mathfrak A}_{[z]}$ such that
\[ N_{[z]} = \cos(t)Z_1 + \sin(t)JZ_2 \]
for some orthonormal vectors $Z_1,Z_2 \in V(A)$ and $0 < t \leq \frac{\pi}{4}$. This implies
\begin{eqnarray*}
N_{[z]} & = & \cos(t)Z_1 + \sin(t)JZ_2, \\
AN_{[z]} & = & \cos(t)Z_1 - \sin(t)JZ_2, \\
\xi_{[z]} & = & \sin(t)Z_2 - \cos(t)JZ_1, \\
A\xi_{[z]} & = & \sin(t)Z_2 + \cos(t)JZ_1,
\end{eqnarray*}
and therefore ${\mathcal Q}_{[z]} = T_{[z]}Q^m \ominus ([Z_1] \oplus [Z_2])$ is strictly contained in ${\mathcal C}_{[z]}$. Moreover, we have
\[
A\xi_{[z]} = B\xi_{[z]}\ {\rm and}\  \rho(\xi_{[z]}) = 0.
\]
We have
\begin{eqnarray*}
g(B\xi_{[z]} + \delta \xi_{[z]},N_{[z]}) & = & 0,\\
g(B\xi_{[z]} + \delta \xi_{[z]},\xi_{[z]}) & = & 0,\\
g(B\xi_{[z]} + \delta\xi_{[z]},B\xi_{[z]} + \delta\xi_{[z]}) & = & \sin^2(2t),
\end{eqnarray*}
where the function $\delta$ denotes ${\delta}=-g({\xi},A{\xi})=-(\sin^2 t - \cos^2 t)=\cos 2t$.
Therefore
\[ U_{[z]} = \frac{1}{\sin(2t)}(B\xi_{[z]} + \delta\xi_{[z]}) \]
is a unit vector in ${\mathcal C}_{[z]}$ and
\[
{\mathcal C}_{[z]} = {\mathcal Q}_{[z]} \oplus [U_{[z]}]\ {\rm (orthogonal\ sum)}.
\]
If $N_{[z]}$ is not ${\mathfrak A}$-principal at $[z]$, then $N$ is not ${\mathfrak A}$-principal in an open neighborhood of $[z]$, and therefore $U$ is a well-defined unit vector field on that open neighborhood. We summarize this in the following

\begin{lemma} \label{notprincipal}
Let $M$ be a real hypersurface in ${Q^*}^m$ whose unit normal $N_{[z]}$ is not ${\mathfrak A}$-principal at $[z]$. Then there exists an open neighborhood of $[z]$ in $M$ and a section $A$ in ${\mathfrak A}$ on that neighborhood consisting of real structures such that
\begin{enumerate}[\rm (i)]
\item $A\xi = B\xi$ and $\rho(\xi) = 0$,
\item $U = (B\xi + \delta\xi)/||B\xi + \delta\xi||$ is a unit vector field tangent to ${\mathcal C}$
\item ${\mathcal C} = {\mathcal Q} \oplus [U]$.
\end{enumerate}
\end{lemma}

\section {Tubes around the totally geodesic ${\mathbb C}H^k \subset {Q^*}^{2k}$}

We assume that $m$ is even, say $m = 2k$. The map
$$
{\mathbb C}H^k \to {Q^*}^{2k}=SO^{o}_{2,2k}/SO_2SO_{2k} \subset {\mathbb C}H^{2k+1} \subset {\mathbb C}_1^{2k+2} ,
$$
is defined by $[z_1,\ldots,z_{k+1}] \mapsto [z_1,\ldots,z_{k+1},iz_1,\ldots,iz_{k+1}]$,
provides an embedding of ${\mathbb C}H^k$ into ${Q^*}^{2k}$ as a totally geodesic complex submanifold in ${\mathbb C}H^{2k+1}$, where
$${Q^*}^{2k}=\{[z_1, {\cdots}, z_{2k+2}]{\in}{\mathbb C}H^{2k+1}{\vert} -z_1^2+z_2^2+{\cdots}z_{k+1}^2-z_{k+2}^2+z_{k+3}^2+{\cdots}+z_{2k+2}^2=0\}$$
can be regarded as the set of negative $2$-planes in indefinite Euclidean space
${\mathbb R}_2^{2k+2}$, that is, a real hyperbolic Grassmannian manifold. Of course, it can be easily checked that the point $[z_1,\ldots,z_{k+1},iz_1,\ldots,iz_{k+1}]$
belongs to ${Q^*}^{2k}$.

Consider the standard embedding of $U_{1,k}$ into $SO^{o}_{2,2k}$ which is determined by the Lie algebra embedding in such a way that
$${\frak u}_{1,k} {\to} {\frak so}_{2,2k},\quad  C+Di {\to} \begin{pmatrix} C & -D \\ D & C \end{pmatrix}\ ,$$
where $C,D{\in}M_{k+1,k+1}({\mathbb R})$ which satisfy respectively ${}^{t}CgC=g$ and ${}^tDgD=g$ for the signature $(1,k)$ of the indefinite Riemannian metric $g$ on ${\mathbb R}_1^{k+1}$ defined by
$g(X,Y)=-x_1y_1+x_2y_2+{\cdots}+x_{k+1}y_{k+1}$ for any $X,Y{\in}{\mathbb R}^{k+1}$.

We define a complex structure $j$ on ${\mathbb C}_1^{2k+2}$ by
$$
j(z_1,\ldots,z_{k+1},z_{k+2},\ldots,z_{2k+2}) = (-z_{k+2},\ldots,-z_{2k+2},z_1,\ldots,z_{k+1}),
$$
where $V(A_{\bar{z}}) =  {\mathbb R}_1^{2k+2} \cap T_{[z]}{Q^*}^{2k}$ is the $(+1)$-eigenspace and $JV(A_{\bar{z}}) = i{\mathbb R}_1^{2k+2} \cap T_{[z]}{Q^*}^{2k}$ is the $(-1)$-eigenspace of $A_{\bar{z}}$.  Geometrically this means that the shape operator $A_{\bar{z}}$ defines a real structure on the complex vector space $T_{[z]}{Q^*}^{2k}$.

Note that $ij = ji$. We can then identify ${\mathbb C}_1^{2k+2}$ with ${\mathbb C}_1^{k+1} \oplus j{\mathbb C}^{k+1}$ and get
$$
T_{[z]}{\mathbb C}H^k = \{X + jiX \mid X \in {\mathbb C}_1^{k+1} \ominus [z]\} = \{X + ijX \mid X \in V(A_{\bar{z}})\}.
$$
Note that the complex structure $j$ on ${\mathbb C}_1^{2k+2}$ corresponds to the complex structure $J$ on $T_{[z]}Q^{2k}$ via the obvious identifications. For the normal space $\nu_{[z]}{\mathbb C}H^k$ of ${\mathbb C}H^k$ at $[z]$ we have
\[
\nu_{[z]}{\mathbb C}H^k = A_{\bar{z}}(T_{[z]}{\mathbb C}H^k) = \{X - ijX \mid X \in V(A_{\bar{z}})\}.
\]
It is easy to see that both the tangent bundle and the normal bundle of ${\mathbb C}H^k$ consist of ${\mathfrak A}$-isotropic singular tangent vectors of ${Q^*}^{2k}$.

We will now calculate the principal curvatures and principal curvature spaces of the tube around ${\mathbb C}H^k$ in ${Q^*}^{2k}$. Let $N$ be a unit normal vector of ${\mathbb C}H^k$ in ${Q^*}^{2k}$ at ${[z]} \in {\mathbb C}H^k$. Since $N$ is ${\mathfrak A}$-isotropic, the four vectors $N,JN,AN,JAN$ are pairwise orthonormal and the normal Jacobi operator ${\bar R}_N$ is given by
\begin{eqnarray*}
{\bar R}_N Z = {\bar R}(Z,N)N  & =   & -Z + g(Z,N)N - 3g(Z,JN)JN  \\
& &  + \, g(Z,AN)AN + g(Z,JAN)JAN.
\end{eqnarray*}
From this, by using that $N$ is $\frak A$-isotropic, ${\bar R}_NN={\bar R}(N,N)N=0$, ${\bar R}_NAN={\bar R}(AN,N)N=0$, ${\bar R}_NJAN=0$, and ${\bar R}_NJN=-4JN$. This implies readily that ${\bar R}_N$ has the three eigenvalues $0,-1$ and $-4$ with corresponding eigenspaces ${\mathbb R}N\oplus [AN]$, $T_{[z]}{Q^*}^{2k} \ominus ([N] \oplus [AN])$ and ${\mathbb R}JN$. Since $[N] \subset \nu_{[z]}{\mathbb C}H^k$ and $[AN] \subset  T_{[z]}{\mathbb C}H^k$, we conclude that both $T_{[z]}{\mathbb C}H^k$ and $\nu_{[z]}{\mathbb C}H^k$ are invariant under ${\bar R}_N$.

To calculate the principal curvatures of the tube around ${\mathbb C}H^k$ we use the Jacobi field method. Let $\gamma$ be the geodesic  in ${Q^*}^{2k}$ with $\gamma(0) = [z]$ and $\dot{\gamma}(0) = N$ and denote by $\gamma^\perp$ the parallel subbundle of $TQ^{2k}$ along $\gamma$ defined by $\gamma^\perp_{\gamma(t)} = T_{[\gamma(t)]}Q^{2k} \ominus {\mathbb R}\dot{\gamma}(t)$. Moreover, define the $\gamma^\perp$-valued tensor field $R^\perp_{\gamma}$ along $\gamma$ by $R^\perp_{\gamma(t)}X = R(X,\dot{\gamma}(t))\dot{\gamma}(t)$.
Now consider the ${\rm End}(\gamma^\perp)$-valued differential equation
\[
Y^{\prime\prime} + R^\perp_{\gamma} \circ Y = 0.
\]
Let $D$ be the unique solution of this differential equation with initial values
\[
D(0) = \begin{pmatrix} I & 0 \\ 0 & 0 \end{pmatrix}\ ,\ D^\prime(0) = \begin{pmatrix} 0 & 0 \\ 0 & I \end{pmatrix},
\]
where the decomposition of the matrices is with respect to
\[
\gamma^\perp_{[z]} = T_{[z]}{\mathbb C}H^k \oplus (\nu_{[z]}{\mathbb C}H^k \ominus {\mathbb R}N)
\]
and $I$ denotes the identity transformation on the corresponding space. Then the shape operator $S(r)$ of the tube around ${\mathbb C}H^k$ with respect to $-\dot{\gamma}(r)$ is given by
\[
S(r) = D^\prime(r)\circ D^{-1}(r).
\]
If we decompose $\gamma^\perp_{[z]}$ further into
\[
\gamma^\perp_{[z]} = (T_{[z]}{\mathbb C}H^k \ominus [AN]) \oplus [AN] \oplus (\nu_{[z]}{\mathbb C}H^k \ominus [N]) \oplus {\mathbb R}JN,
\]
we get by explicit computation that
\[
S(r) = \begin{pmatrix}
\tanh(r) & 0 & 0 & 0 \\
0 & 0 & 0 & 0 \\
0 & 0 & \coth(r) & 0 \\
0 & 0 & 0 & 2\coth(2r)
 \end{pmatrix}
\]
with respect to that decomposition. Here let us check that $SJN=2\coth(2r)JN$ for $M{\subset}{Q^*}^{2k}$. Since ${\bar R}_NJN=-4JN$, we have $Y^{''}-4Y=0$ for a geodesic $\gamma$ such that ${\gamma}(0)=[z]$ and ${\gamma}'(0)=N$.
The the solution vector field $Y(r)$ of the Jacobi equation becomes
$$Y(r)=(c_1\cosh(2r)+c_2\sinh(2r))E_X(r).$$
By the initial condition $0=Y(0)=c_1E_X(0)=c_1X$ and $X=Y'(0)=2c_2E_X(0)=2c_2X$, we know that the solution vector field is given by $Y(r)=D(r)E_X(r)=\frac{1}{2}\sinh(2r)E_X(r)$. From this, together with the definition of the shape operator, it follows that
\begin{equation*}
\begin{split}
\frac{1}{2}\sinh(2r)S(r)E_X(r)=&S(r)Y(r)=D'(r)D^{-1}(r)Y\\
=&D'(r)E_X(r)=\cosh(2r)E_X(r).
\end{split}
\end{equation*}
This implies that $S(r)E_X(r)=2\coth(2r)E_X(r)$, which means $S(r)JN=2\coth(2r)JN$.  By using the similar method we can calulate the other principal curvatures.
Therefore the tube around ${\mathbb C}H^k$ has four distinct constant principal curvatures $\tanh(r)$, $0$, $\coth(r)$ and $2\coth(2r)$ (unless $m=2$ in which case there are only two distinct constant principal curvatures $0$ and $2\coth(2r)$). The corresponding principal curvature spaces are  $T_{[z]}{\mathbb C}H^k \ominus [AN]$, $[AN]$, $\nu_{[z]}{\mathbb C}H^k \ominus [N]$ and ${\mathbb R}JN$ respectively, where we identify the subspaces obtained by parallel translation along $\gamma$ from $[z]$ to $\gamma(r)$.

Note that the parallel translate of $[AN]$ corresponds to ${\mathcal C} \ominus {\mathcal Q}$, the parallel translate of $[N]$ corresponds to ${\mathbb C}\nu M$, and the parallel translate of ${\mathbb R}JN$ corresponds to ${\mathcal F}$. Moreover, we have $A(T_{[z]}{\mathbb C}H^k \ominus [AN]) = \nu_{[z]}{\mathbb C}H^k \ominus [N]$.

When $M$ becomes an open part of a horosphere in ${Q^*}^{2k}$ whose center at infinity in the equivalence class of an $\mathfrak A$-isotropic geodesic in ${Q^*}^{2k}$, by using the results in Suh \cite{S5}, \cite{S6} and \cite{S7} we can calculate that it has three distinct constant prinicipal curvatures $1,0,1$ and $2$ corresponding to the same principal curvature spaces mentioned above.

Since $JN$ is a principal curvature vector, we conclude that every tube around ${\mathbb C}H^k$ is a Hopf hypersurface. We also see that all principal curvature spaces orthogonal to ${\mathbb R}JN$ are $J$-invariant. Thus, if $\phi$ denotes the structure tensor field on the tube which is induced by $J$, we get $S\phi = \phi S$. Since the complex structure on $Q^m$ is parallel, we have
\[
g(\nabla_X\xi,Y) + g(X,\nabla_Y\xi) = g((S\phi - \phi S)X,Y)
\]
for all $X,Y \in TM$. As $\xi$ is a Killing vector field if and only
if $\nabla\xi$ is a skew-symmetric tensor field, we see that the Reeb flow on $M$ is isometric if and only if $S\phi = \phi S$.

We summarize the previous discussion in the following proposition.
\par
\vskip 6pt

\begin{proposition}
Let $M$ be the tube around the totally geodesic ${\mathbb C}H^k$ in ${Q^*}^{2k}$, $k \geq 2$, or the horosphere in ${Q^*}^{2k}$ whose center at infinity is in the equivalent class of an $\frak A$-isotropic singular geodesic in ${Q^*}^{2k}$. Then the following statements hold:
\begin{enumerate}[\rm (i)]
\item {$M$ is a Hopf hypersurface.}
\item {The tangent bundle $TM$ and the normal bundle $\nu M$ of $M$ consist of ${\mathfrak A}$-isotropic singular tangent vectors of ${Q^*}^{2k}$.}
\item {$M$ has four(or three) distinct constant principal curvatures. Their values and corresponding principal curvature spaces and multiplicities are given in the following table:
\begin{table}[ht]
\caption{Principal curvatures of $M$}
{\begin{tabular}{|c|c|c|}
\hline
{\rm principal curvature} & {\rm eigenspace} & {\rm multiplicity} \\
\hline
$2\coth(2r), 2$ & ${\mathcal F}$ & $1$ \\
$0$ & ${\mathcal C} \ominus {\mathcal Q}$ & $2$ \\
$\tanh(r), 1$ & $T{\mathbb C}P^k \ominus ({\mathcal C} \ominus {\mathcal Q})$ & $2k-2$\\
$\coth(r), 1$ & $\nu{\mathbb C}P^k \ominus {\mathbb C}\nu M$ & $2k-2$ \\
\hline
\end{tabular}}
\end{table}
}
\noindent The real structure $A$ determined by the ${\mathfrak A}$-isotropic unit normal vector at $[z]$ maps $T_{[z]}{\mathbb C}H^k \ominus ({\mathcal C}_{[z]} \ominus {\mathcal Q}_{[z]})$ onto $\nu_{[z]}{\mathbb C}H^k \ominus {\mathbb C}\nu_{[z]} M$, and vice versa.
\item {The shape operator $S$ of $M$ and the structure tensor field $\phi$ of $M$ commute with each other, that is, $S\phi = \phi S$.}
\item {The Reeb flow on $M$ is an isometric flow.}
\end{enumerate}
\end{proposition}

\section {The Codazzi equation and some consequences}

\medskip
From the explicit expression of the Riemannian curvature tensor of the complex quadric $Q^m$ we can easily derive the Codazzi equation for a real hypersurface $M \subset {Q^*}^m$:
\begin{eqnarray*}
& & g((\nabla_XS)Y - (\nabla_YS)X,Z) \\
& = &- \eta(X)g(\phi Y,Z) + \eta(Y) g(\phi X,Z) + 2\eta(Z) g(\phi X,Y) \\
 & & -\, \rho(X)g(BY,Z) + \rho(Y)g(BX,Z)\\
& & +\, \eta(BX)g(BY,\phi Z) + \eta(BX)\rho(Y)\eta(Z)  \\
& & -\, \eta(BY)g(BX,\phi Z) - \eta(BY)\rho(X)\eta(Z)  .
\end{eqnarray*}
We now assume that $M$ is a Hopf hypersurface. Then we have
\[
S\xi = \alpha \xi
\]
with the smooth function $\alpha = g(S\xi,\xi)$ on $M$.
Inserting $Z = \xi$ into the Codazzi equation leads to
\[
g((\nabla_XS)Y - (\nabla_YS)X,\xi) =   2 g(\phi X,Y)  - 2\rho(X)\eta(BY) + 2\rho(Y)\eta(BX).
\]
On the other hand, we have
\begin{eqnarray*}
 & & g((\nabla_XS)Y - (\nabla_YS)X,\xi) \\
& = & g((\nabla_XS)\xi,Y) - g((\nabla_YS)\xi,X) \\
& = & d\alpha(X)\eta(Y) - d\alpha(Y)\eta(X) + \alpha g((S\phi + \phi S)X,Y) - 2g(S \phi SX,Y).
\end{eqnarray*}
Comparing the previous two equations and putting $X = \xi$ yields
$$
d\alpha(Y)  =  d\alpha(\xi)\eta(Y)  + 2\delta\rho(Y).
$$
Reinserting this into the previous equation yields
\begin{eqnarray*}
 & & g((\nabla_XS)Y - (\nabla_YS)X,\xi) \\
& = & -2\delta\eta(X)\rho(Y) + 2\delta \rho(X)\eta(Y) \\
& & +\, \alpha g((\phi S + S\phi)X,Y) - 2g(S \phi SX,Y) .
\end{eqnarray*}
Altogether this implies
\begin{eqnarray*}
0 & = & 2g(S \phi SX,Y) - \alpha g((\phi S + S\phi)X,Y) + 2 g(\phi X,Y) \\
& & - 2\delta \rho(X)\eta(Y) - 2\rho(X)\eta(BY) + 2\rho(Y)\eta(BX) + 2\delta\eta(X)\rho(Y) \\
& = & g((2S \phi S - \alpha(\phi S + S\phi) + 2\phi) X ,Y)\\
& & - 2\rho(X)\eta(BY+\delta Y) + 2\rho(Y)\eta(BX+\delta X) \\
& = & g((2S \phi S - \alpha(\phi S + S\phi) + 2\phi) X ,Y)\\
& & - 2\rho(X)g(Y,B\xi+\delta \xi) + 2g(X,B\xi+\delta \xi)\rho(Y).
\end{eqnarray*}
If $AN = N$ we have $\rho = 0$, otherwise we can use Lemma \ref{notprincipal} to calculate
$\rho(Y) = g(Y,AN) = g(Y,AJ\xi) = -g(Y,JA\xi) = -g(Y,JB\xi) = -g(Y,\phi B\xi)$.
Thus we have proved

\begin{lemma}\label{Cod1}
Let $M$ be a Hopf hypersurface in ${Q^*}^m$, $m \geq 3$. Then we have
\[
(2S \phi S - \alpha(\phi S + S\phi) + 2\phi) X =
 2\rho(X)(B\xi+\delta \xi) + 2g(X,B\xi+\delta \xi)\phi B\xi .
\]
\end{lemma}

If the unit normal vector field $N$ is ${\mathfrak A}$-principal, we can choose a real structure $A \in {\mathfrak A}$ such that $AN = N$. Then we have $\rho = 0$ and $\phi B\xi = -\phi \xi = 0$, and therefore
\[
2S \phi S - \alpha(\phi S + S\phi) = -2\phi.
\]
If $N$ is not ${\mathfrak A}$-principal, we can choose a real structure $A \in {\mathfrak A}$ as in Lemma \ref{notprincipal} and get
\begin{eqnarray*}
& & \rho(X)(B\xi+\delta \xi) + g(X,B\xi+\delta \xi)\phi B\xi \\
& = & -g(X,\phi(B\xi + \delta\xi))(B\xi+\delta \xi) + g(X,B\xi+\delta \xi)\phi (B\xi + \delta\xi)\\
& = & ||B\xi + \delta\xi||^2( g(X,U)\phi U -g(X,\phi U)U ) \\
& = & \sin^2(2t)( g(X,U)\phi U -g(X,\phi U)U ),
\end{eqnarray*}
which is equal to $0$ on ${\mathcal Q}$ and equal to $\sin^2(2t)\phi X$ on ${\mathcal C} \ominus {\mathcal Q}$. Altogether we have proved:

\begin{lemma}\label{Cod2}
Let $M$ be a Hopf hypersurface in ${Q^*}^m$, $m \geq 3$. Then the tensor field
\[ 2S\phi S - \alpha (\phi S + S\phi) \]
leaves ${\mathcal Q}$ and ${\mathcal C} \ominus {\mathcal Q}$ invariant and we have
\[ 2S\phi S - \alpha (\phi S + S\phi) = -2\phi \ {\rm on}\  {\mathcal Q} \]
and
\[ 2S\phi S - \alpha (\phi S + S\phi) = -2\delta^2\phi \ {\rm on}\  {\mathcal C} \ominus {\mathcal Q}, \]
where ${\delta}=\cos 2t$ as in section 3.
\end{lemma}

We will now prove that the principal curvature $\alpha$ of a Hopf hypersurface is constant if the normal vectors are ${\mathfrak A}$-isotropic. Assume that the unit normal vector field $N$ is ${\mathfrak A}$-isotropic everywhere. Then we have $\delta = 0$ and we get
\[
Y\alpha  =  d\alpha(\xi)\eta(Y)
\]
for all $Y \in TM$.
Since ${\rm grad}^M \alpha =d\alpha(\xi)\xi$, we can compute the Hessian ${\rm hess}^M \alpha$ by
\begin{eqnarray*}
({\rm hess}^M \alpha)(X,Y) & = & g(\nabla_X{\rm grad}^M \alpha , Y) \\
& = & d(d\alpha(\xi))(X)\eta(Y) +d\alpha(\xi) g(\phi SX,Y).
\end{eqnarray*}
As ${\rm hess}^M \alpha$ is a symmetric bilinear form, the previous equation implies
\[
d\alpha(\xi)g((S\phi + \phi S)X,Y) = 0
\]
for all vector fields $X,Y$ on $M$ which are tangential to ${\mathcal C}$.

Now let us assume that $S{\phi}+{\phi}S=0$. For every principal curvature vector $X \in {\mathcal C}$ such that $SX={\lambda}X$ this implies $S{\phi}X=-{\phi}SX=-{\lambda}{\phi}X$. We assume $||X|| = 1$ and put $Y = \phi X$. Using the normal vector field $N$ is $\frak A$-isotropic, that is ${\delta}=0$ in Lemma \ref{Cod1} we know that
$$-{\lambda}^2{\phi}X+{\phi}X={\rho}(X)B{\xi}+g(X,B{\xi}){\phi}B{\xi}.$$
From this, using $g(X,B{\xi})=g(X,A{\xi})=-g({\phi}X,AN)=-{\rho}({\phi}X)$, and together with the fact that ${\lambda}=0$ from the above equation and the commuting shape operator $$-{\lambda}{\phi}X=-{\phi}SX=S{\phi}X={\phi}SX={\lambda}{\phi}X,$$
we get ${\lambda}=0$. This implies that
\begin{eqnarray*}
1 & = & -\lambda^2 + 1 = \rho(X)\eta(B\phi X) - \rho(\phi X)\eta(BX)\\
& =&  g(X,AN)^2 + g(X,A\xi)^2 =  ||X_{{\mathcal C} \ominus {\mathcal Q}}||^2 \leq 1,
\end{eqnarray*}
where $X_{{\mathcal C} \ominus {\mathcal Q}}$ denotes the orthogonal projection of $X$ onto ${\mathcal C} \ominus {\mathcal Q}$. This
implies $||X_{{\mathcal C} \ominus {\mathcal Q}}||^2 = 1$
for all principal curvature vectors $X \in {\mathcal C}$ with $||X|| = 1$. This is only possible if ${\mathcal C} = {\mathcal C} \ominus {\mathcal Q}$, or equivalently, if ${\mathcal Q} = 0$. Since $m \geq 3$ this is not possible. Hence we must have $S{\phi}+{\phi}S \neq 0$ everywhere and therefore $d\alpha(\xi) = 0$, which implies ${\rm grad}^M \alpha = 0$. Since $M$ is connected this implies that $\alpha$ is constant. Thus we have proved:

\medskip

\begin{lemma} \label{alphaconstant}
Let $M$ be a real hypersurface in ${Q^*}^m$, $m \geq 3$, with isometric Reeb flow and ${\mathfrak A}$-isotropic normal vector
field $N$ everywhere. Then $\alpha$ is constant.
\end{lemma}

\section{Proof of Theorem \ref{mainresult}}

Let $M$ be a real hypersurface in ${Q^*}^m$, $m \geq 3$, with isometric Reeb flow. As we have seen above, this geometric condition is equivalent to the algebraic condition  $S\phi = \phi S$. Applying this equation to $\xi$ gives $0 = S\phi \xi = \phi S\xi$, which implies that $S\xi = \alpha \xi$ with $\alpha = g(S\xi,\xi)$. Therefore any real hypersurface in ${Q^*}^m$ with isometric Reeb flow is a Hopf hypersurface.

Differentiating the equation $S\phi - \phi S = 0$ gives
\begin{eqnarray*}
0 & = & (\nabla_XS)\phi Y + S(\nabla_X\phi)Y - (\nabla_X\phi)SY - \phi(\nabla_XS)Y \\
& = & (\nabla_XS)\phi Y + S(\eta(Y)SX-g(SX,Y)\xi) \\
& & -\, (\eta(SY)SX-g(SX,SY)\xi) - \phi(\nabla_XS)Y \\
& = & (\nabla_XS)\phi Y + \eta(Y)S^2X-\alpha g(SX,Y)\xi \\
& & -\, \eta(SY)SX + g(SX,SY)\xi - \phi(\nabla_XS)Y.
\end{eqnarray*}
If we define
\[
\Theta(X,Y,Z) = g((\nabla_XS)Y,\phi Z) + g((\nabla_XS)Z,\phi Y),
\]
the previous equation implies
\begin{eqnarray*}
\Theta(X,Y,Z) & = & \alpha\eta(Z)g(SX,Y) - \eta(Z)g(SX,SY) \\
& & +\, \eta(SY)g(SX,Z) - \eta(Y)g(S^2X,Z).
\end{eqnarray*}
Evaluating $\Theta(X,Y,Z) + \Theta(Y,Z,X) - \Theta(Z,X,Y)$  leads to
\begin{eqnarray*}
2g((\nabla_XS)Y,\phi Z) & = & \Phi(X,Y,Z) - \Phi(Y,Z,X) + \Phi(Z,X,Y) \\
& & + 2\alpha\eta(Z)g(SX,Y) - 2\eta(Z)g(S^2X,Y),
\end{eqnarray*}
where
\[
\Phi(X,Y,Z) = g((\nabla_XS)Y - (\nabla_YS)X,\phi Z).
\]
The three $\Psi$-terms can be evaluated using the Codazzi equation, which leads to
\begin{eqnarray*}
2g((\nabla_XS)Y,\phi Z)
& = &  \rho(X)\Big\{g(AY,\phi Z) - g(AZ,\phi Y)\Big\} \\
& & + \eta(BX)\Big\{g(JAY,\phi Z) - g(JAZ,\phi Y)\Big\}\\
& & - \rho(Y)\Big\{g(AX,\phi Z) + g(AZ,\phi X)\Big\} \\
& & -  \eta(BY)\Big\{g(JAX,\phi Z) + g(JAZ,\phi X)\Big\}\\
& & + \rho(Z)\Big\{g(AX,\phi Y) + g(AY,\phi X)\Big\} \\
& & + \eta(BZ)(g(JAX,\phi Y) + g(JAY,\phi X))\\
& & + 2\eta(Z)g(\phi X,\phi Y) - 2\eta(Y)g(\phi X,\phi Z) \\
& & + 2\alpha\eta(Z)g(SX,Y) - 2\eta(Z)g(S^2X,Y).
\end{eqnarray*}
Replacing $\phi Z$ by $JZ - \eta(Z)N$, and similarly for $X$ and $Y$, one can easily calculate that
\begin{eqnarray*}
{\eta}(BX)\{g(JAY,\phi Z) - g(JAZ,\phi Y)\} & = &{\eta}(BX)\Big\{ \eta(Y)\eta(BZ) - \eta(Z)\eta(BY)\Big\}, \\
-{\eta}(BY)\{g(JAX,\phi Z) + g(JAZ,\phi X)\} & = & -\Big\{2g(BX,Z) - \eta(X)\eta(BZ) \\
& & - \eta(Z)\eta(BX)\Big\}{\eta}(BY),\\
{\eta}(BZ)\{g(JAX,\phi Y) + g(JAY,\phi X) & = &\Big\{2g(BX,Y) - \eta(X)\eta(BY) \\
& & - \eta(Y)\eta(BX)\Big\}{\eta}(BZ).
\end{eqnarray*}
Inserting this into the previous equation gives
\begin{eqnarray*}
2g((\nabla_XS)Y,\phi Z)
& = &  - \rho(X)\Big\{g(BY,\phi Z) - g(BZ,\phi Y)\Big\} \\
& & + \rho(Y)\Big\{g(BX,\phi Z) + g(BZ,\phi X)\Big\} \\
& & - \rho(Z)\Big\{g(BX,\phi Y) + g(BY,\phi X)\Big\} \\
& & + 2\eta(BY)g(BX,Z) - 2\eta(BZ)g(BX,Y)  \\
& & - 2\eta(Z)g(\phi X,\phi Y) + 2\eta(Y)g(\phi X,\phi Z) \\
& & + 2\alpha\eta(Z)g(SX,Y) - 2\eta(Z)g(S^2X,Y).
\end{eqnarray*}
Since
\begin{eqnarray*}
{\rho}(X)\{g(BY,\phi Z) - g(BZ,\phi Y)\}  & =  &{\rho}(X)\{ \eta(Y)\rho(Z) - \eta(Z)\rho(Y)\}, \\
{\rho}(Z)\{g(BX,\phi Y) + g(BY,\phi X)\}  & =  &\{ 2g(BX,\phi Y) + \eta(Y)\rho(X) \\
& & - \eta(X)\rho(Y)\}{\rho}(Z), \\
-{\rho}(Y)\{g(BX,\phi Z) + g(BZ,\phi X)\}  & =  &-\{ 2g(BX,\phi Z) + \eta(Z)\rho(X) \\
& & - \eta(X)\rho(Z)\}{\rho}(Y),
\end{eqnarray*}
we get
\begin{eqnarray*}
g((\nabla_XS)Y,\phi Z)
& = &  -\rho(X)\eta(Y)\rho(Z) + \rho(X)\rho(Y)\eta(Z) \\
& & + \rho(Y)g(BX,\phi Z) - \rho(Z)g(BX,\phi Y) \\
& &  + \eta(BY)g(BX,Z) - \eta(BZ)g(BX,Y)  \\
& & - \eta(Z)g(\phi X,\phi Y) + \eta(Y)g(\phi X,\phi Z) \\
& & + \alpha\eta(Z)g(SX,Y) - \eta(Z)g(S^2X,Y).
\end{eqnarray*}
Replacing $Z$ by $\phi Z$ and using $\phi^2 Z = - Z + \eta(Z)\xi$ gives
\begin{eqnarray*}
& & -g((\nabla_XS)Y,Z) + \eta(Z)g((\nabla_XS)Y,\xi) \\
& = & - \rho(X)\eta(Y)\rho(\phi Z)  - \rho(Y)g(BX,Z) \\
& & + \rho(Y)\eta(Z)g(BX,\xi) - \rho(\phi Z)g(BX,\phi Y) \\
& &  + \eta(BY)g(BX,\phi Z) - \eta(B\phi Z)g(BX,Y) - \eta(Y)g(\phi X,Z).
\end{eqnarray*}
Since
\[
g((\nabla_XS)Y,\xi) = d\alpha(X)\eta(Y) + \alpha g(S \phi X,Y) - g(S^2 \phi X,Y),
\]
this implies
\begin{eqnarray*}
g((\nabla_XS)Y,Z)
& = &  d\alpha(X)\eta(Y)\eta(Z) + \eta(Z) g((\alpha S \phi - S^2 \phi) X,Y) \\
& &  + \rho(X)\eta(Y)\rho(\phi Z)  + \rho(Y)g(BX,Z) \\
& & - \rho(Y)\eta(Z)g(BX,\xi) + \rho(\phi Z)g(BX,\phi Y) \\
& &  - \eta(BY)g(BX,\phi Z) + \eta(B\phi Z)g(BX,Y) \\
& & + \eta(Y)g(\phi X,Z).
\end{eqnarray*}
From this we get an explicit expression for the covariant derivative of the shape operator,
\begin{eqnarray*}
(\nabla_XS)Y  & = &
\Big\{d\alpha(X)\eta(Y) + g((\alpha S \phi - S^2 \phi) X,Y)  - \delta\eta(Y)\rho(X)  \\
& & \quad - \delta g(BX,\phi Y) - \eta(BX)\rho(Y)\Big\} \xi\\
& &  - \{\eta(Y)\rho(X) +g(BX,\phi Y)\} B\xi - g(BX,Y)\phi B\xi  \\
& & + \rho(Y)BX + \eta(Y)\phi X + \eta(BY)\phi BX  .
\end{eqnarray*}
Putting $Y = \xi$ and $X \in {\mathcal C}$ then leads to
\begin{eqnarray*}
\alpha S \phi X - S^2 \phi  X & = & - \delta \rho(X) \xi - \rho(X)  B\xi - \eta(BX)\phi B\xi + \phi X - \delta\phi BX  .
\end{eqnarray*}
On the other hand, from Lemma \ref{Cod1} we get
\begin{eqnarray*}
\alpha S \phi X - S^2 \phi  X & = &  \delta \rho(X) \xi + \rho(X)  B\xi + \eta(BX)\phi B\xi - \phi X   .
\end{eqnarray*}
Comparing the previous two equations leads to
\[
\delta\phi BX = 0
\]
for all $X \in {\mathcal C}$. Let us first assume that $\delta \neq 0$. Then we have $\phi BX = 0$ for all $X \in {\mathcal C}$, which implies $BX = \eta(BX)\xi$ for all $X \in {\mathcal C}$, and therefore
\[
AX = BX + \rho(X)N = \eta(BX)\xi + \rho(X)N
\]
for all $X \in {\mathcal C}$. This implies $A({\mathcal C}) \subset [N]$, which gives a contradiction since $A$ is an isomorphism everywhere and the rank of ${\mathcal C}$ is equal to $2(m-1)$ and $m \geq 3$. Therefore we must have $\delta = 0$, which means that $N$ is ${\mathfrak A}$-isotropic. We thus have proved the following proposition.

\begin{proposition} \label{Aisotropic}
Let $M$ be a real hypersurface in ${Q^*}^m$, $m \geq 3$, with isometric Reeb flow. Then the normal vector
field $N$ is ${\mathfrak A}$-isotropic everywhere.
\end{proposition}

From Proposition \ref{Aisotropic} and Lemma \ref{alphaconstant} we conclude that the principal curvature function $\alpha$ is constant.
Since the normal vector field $N$ is ${\mathfrak A}$-isotropic, the vector fields $N$, $\xi$, $B\xi$ and $\phi B\xi$ are pairwise orthonormal. This implies that
\[
{\mathcal C} \ominus {\mathcal Q} = [B\xi].
\]
From Lemma \ref{Cod2} we know that the tensor field $2S\phi S - \alpha (\phi S + S\phi) = 2\phi (S^2 - \alpha S)$
leaves ${\mathcal Q}$ and ${\mathcal C} \ominus {\mathcal Q}$ invariant and
\[ \phi (S^2 - \alpha S) = \phi \   {\rm on}\  {\mathcal Q} \ \ {\rm and} \ \
\phi (S^2 - \alpha S) = 0 \ {\rm on}\  {\mathcal C} \ominus {\mathcal Q}. \]
Since $\phi$ is an isomorphism of ${\mathcal Q}$ and of ${\mathcal C} \ominus {\mathcal Q}$ this implies
\[ S^2 - \alpha S = -I \ {\rm on}\  {\mathcal Q} \ \ {\rm and}\ \
S^2 - \alpha S = 0 \ {\rm on}\  {\mathcal C} \ominus {\mathcal Q}. \]
As $M$ is a Hopf hypersurface we have $S({\mathcal C}) \subset {\mathcal C}$. Let $X \in {\mathcal C}$ be a principal curvature vector of $M$ with corresponding principal curvature $\lambda$, that is, $SX = \lambda X$. We decompose $X$ into $X = Y + Z$ with $Y \in {\mathcal Q}$ and $Z \in {\mathcal C} \ominus {\mathcal Q}$. Then we get
\[
(\lambda^2 -\alpha\lambda)Y + (\lambda^2 -\alpha\lambda)Z = (\lambda^2 -\alpha\lambda)X = (S^2 - \alpha S)X = -Y.
\]
If $\lambda^2 - \alpha\lambda = 0$ we must have $Y = 0$ and therefore $X \in {\mathcal C} \ominus {\mathcal Q}$. If $\lambda^2 - \alpha\lambda \neq 0$ we must have $Z = 0$ and therefore $X \in {\mathcal Q}$. Altogether this implies

\begin{proposition} \label{Qinvariance}
Let $M$ be a real hypersurface in ${Q^*}^m$, $m \geq 3$, with isometric Reeb flow. Then the distributions ${\mathcal Q}$ and ${\mathcal C} \ominus {\mathcal Q} = [B\xi]$ are invariant under the shape operator $S$ of $M$.
\end{proposition}

Assume that $SX = \lambda X$ with $X \in {\mathcal C}$. From Lemma \ref{Cod2} we get $\lambda^2 - \alpha\lambda + 1 = 0$ if $X \in {\mathcal Q}$ and $\lambda^2 - \alpha\lambda = 0$ if $X \in {\mathcal C} \ominus {\mathcal Q}$. Recall that $\alpha$ is constant. We put $\alpha = 2\coth(2r)$ and define $T_\alpha = {\mathbb R}\xi = {\mathcal F}$. Then the solutions of $x^2 -\alpha x + 1 = 0$ are $\lambda = \coth(r)$ and $\mu = \tanh(r)$. We denote by $T_\lambda$ and $T_\mu$ the subbundles of ${\mathcal Q}$ consisting of the corresponding principal curvature vectors. The rank of ${\mathcal C} \ominus {\mathcal Q}$ is equal to $2$ and ${\mathcal C} \ominus {\mathcal Q}$ is both $S$- and $\phi$-invariant. Therefore, since $S\phi = \phi S$, there is exactly one principal curvature $\kappa$ on ${\mathcal C} \ominus {\mathcal Q}$ which is equal to either $\kappa = 0$ or $\kappa = \alpha$, because the distributions $\mathcal Q$ and ${\mathcal C} \ominus {\mathcal Q}$ are both $S$-invariant and $\phi$-invariant. We define $T_\kappa = {\mathcal C} \ominus {\mathcal Q}$. Note that, since $S\phi = \phi S$, we have $JT_\rho = \phi T_\rho = T_\rho$ for $\rho \in \{\lambda,\mu,\kappa\}$.

\medskip
According to Lemma \ref{Cod1} we have
\begin{eqnarray*}
(\alpha S\phi - S^2 \phi)X & = &
\phi X  - \rho(X)B\xi - g(X,B\xi)\phi B\xi \\
& = &   \phi X  - \eta(B\phi X)B\xi - \eta(BX)\phi B\xi.
\end{eqnarray*}
Inserting this into the above expression for the covariant derivative of $S$ and replacing $\rho(X)$ and $\rho(Y)$ by $\eta(B\phi X)$ and $\eta(B \phi Y)$ respectively leads to
\begin{eqnarray*}
(\nabla_XS)Y  & = &
\{g((\alpha S \phi - S^2 \phi) X,Y)  - \eta(BX)\eta(B\phi Y)\} \xi\\
& &  - \{\eta(Y)\eta(B \phi X) +g(BX,\phi Y)\} B\xi - g(BX,Y)\phi B\xi  \\
& & + \eta(B \phi Y)BX + \eta(Y)\phi X + \eta(BY)\phi BX  \\
& = & \{-\eta(B\phi X)\eta(BY) + g(\phi X,Y)\} \xi\\
& &  - \{\eta(B\phi X)\eta(Y)-g(\phi BX,Y)\} B\xi - g(BX,Y)\phi B\xi  \\
& & + \eta(B\phi Y)BX + \eta(Y)\phi X + \eta(BY)\phi BX  .
\end{eqnarray*}
This implies for any $X,Y{\in}{\mathcal Q}$
\[ -g(BX,Y) = g((\nabla_XS)Y,\phi B\xi).\]
Let $X$ and $Y$ be sections in $T_\lambda$. Then the previous equation implies
\begin{eqnarray*}
-g(BX,Y) & = & (\lambda - \kappa)g(\nabla_XY,\phi B\xi) = (\kappa - \lambda)g(\nabla_X(\phi B\xi),Y) \\
& = & (\kappa - \lambda)g(\phi\nabla_X(B\xi),Y).
\end{eqnarray*}
Since $B\xi = A\xi$, the Gauss formula for $M$ in ${Q^*}^m$ gives
\[
\nabla_X(B\xi) = \nabla_X(A\xi) + g(SX,B{\xi})N = \bar\nabla_X(A\xi),
\]
where $\bar\nabla$ is the Riemannian connection of ${Q^*}^m$. Since ${\mathfrak A}$ is a parallel subbundle of ${\rm End}(T{Q^*}^m)$, by Smith \cite{S}, ${\bar\nabla}_XA=q(X)A$ for a certain $1$-form $q$ on $M$, so there exists $A^\prime \in {\mathfrak A}$ such that $\bar\nabla_XA = A^\prime$, and we get
\[
\bar\nabla_X(A\xi) = A\bar\nabla_X\xi + A^\prime\xi = A\nabla_X\xi + A^\prime\xi = A\phi SX + A^\prime \xi = \lambda A\phi X + A^\prime\xi.
\]
Since $A^\prime$ differs from $A$ by a complex scalar we have $A^\prime \xi =q(X)A{\xi}\in [A\xi] = [B\xi]$ and thus $g(\phi A^\prime \xi, Y) = 0$.
Altogether we therefore get
\begin{eqnarray*}
-g(BX,Y) & = & (\kappa - \lambda)\lambda g(\phi A\phi X,Y) = (\kappa - \lambda)\lambda g(AX,Y) \\
& = & (\kappa - \lambda)\lambda g(BX,Y).
\end{eqnarray*}
Recall that $\kappa \in \{0,\alpha\}$. When $\kappa = 0$ we get $\lambda^2 = 1$ or $g(BX,Y)=0$ for any $X,Y{\in}T_{\lambda}$, ${\lambda}=\coth(r)$, so it follows that ${\lambda}={\mu}=1$, ${\alpha}=2$ or $g(BX,Y)=0$.
\par
When ${\kappa}={\alpha}$, let us suppose that $B(X,Y){\not =}0$. Then we get $(\alpha - \lambda)\lambda = -1$, which implies that $-1={\lambda}{\mu}=1$ . This gives us a contradiction. Therefore we must have $B(X,Y)=0$ for any $X,Y{\in}T_{\lambda}$.  The same argument can be repeated for $T_\mu$. Since ${\mathcal Q} = T_\lambda \oplus T_\mu$ and $B({\mathcal Q}) = {\mathcal Q}$ we conclude
\[
B(T_\lambda) = T_\mu\ \ {\rm and}\ \ B(T_\mu) = T_\lambda.
\]
Since $B = A$ on $T_\lambda$ and $T_\mu$ we can replace $B$ by $A$ here. As both $T_\lambda$ and $T_\mu$ are complex we see that ${\mathcal Q}$ and hence ${Q^*}^m$ must have even complex dimension. We summarize this in

\begin{proposition} \label{Qinvariance}
Let $M$ be a real hypersurface in ${Q^*}^m$, $m \geq 3$, with isometric Reeb flow. Then ${\alpha}=2$, ${\lambda}={\mu}=1$ and ${\kappa}=0$ or $g(BX,Y)=0$ for any $X,Y{\in}T_{\lambda}$, the latter case occurs only for $m = 2k$, that is, the real structure $A$ maps $T_\lambda$ onto $T_\mu$, and vice versa.
\end{proposition}

Let ${\alpha}=2$, ${\lambda}={\mu}=1$, and ${\kappa}=0$. Then the corresponding Jacobi vector field $Y_X$ is resepectively given as

\begin{equation*}
Y_X(r)=\begin{cases}
\exp(-2r)E_X(r) & \text{if $X \in T_\alpha$, ${\alpha}=2$, }\\
\exp(-r)E_X(r) & \text{if $X\in T_\rho$, $\rho\in\{\lambda ,  \mu \}$, $\lambda = \mu =1$, }\\
E_X(r) & \text{if $X\in T_\kappa$, ${\kappa}=0$},
\end{cases}
\end{equation*}

The corresponding Jacobi field is a linear combination of the three types of the Jacobi fields $Y_X$ listed above, and hence its length remains bounded when $r{\to}{\infty}$. This means that all that Jacobi vector fields are bounded for $r{\to}{\infty}$. This shows that all normal geodesics ${\gamma}_{[z]}$ of $M$ in ${Q^*}^m$ are asymtotic to each other and hence determine a singular point $z{\in}{Q^*}^{m}({\infty})$ at infinty(see \cite{E}, \cite{S5} and \cite{S6}). Accordingly, $M$ is an integral manifold of the distribution on the asymtotic class $z$. This distribution is integrable and the maximal leaves are the horospheres in ${Q^*}^m$ whose center at infinity is $z$. Such a uniqueness of integral manifolds of the integrable distributions finally implies that
 $M$ becomes a horosphere in ${Q^*}^m$ whose center is the singular point $z$ at infinity. Moreover, it is given by an equivalence class of assymtotic geodesics whose tangent vectors are $\frak A$-isotropic. Here the meaning of $\frak A$-isotropic means that there exist two orthonormal vector fields $X$ and $Y$ in $V(A)$ such that $AN=\frac{1}{\sqrt 2}X+\frac{1}{\sqrt 2}JY$. Thus the center at infinity of the horosphere is a singular point of type $\mathfrak A$-isotropic.

\par
\vskip 6pt

For each point $[z] \in M$ we denote by $\gamma_{[z]}$ the geodesic in ${Q^*}^{2k}$ with $\gamma_{[z]}(0)=[z]$ and $\dot{\gamma}_{[z]}(0)=N_{[z]}$ and by $F$ the smooth map
\[
F: M \longrightarrow {Q^*}^{2k}, [z] \longrightarrow \gamma_{[z]}(r).
\]
Geometrically, $F$ is the displacement of $M$ at distance $r$ in the direction of the normal vector field $N$. For each point $[z] \in M$ the differential $d_{[z]}F$ of $F$ at $[z]$ can be computed by using Jacobi vector fields as
\[
d_{[z]}F(X)=Z_X(r),
\]
where $Z_X$ is the Jacobi vector field along $\gamma_{[z]}$ with initial values $Z_X(0)=X$ and $Z_X'(0)=-SX$. Using the fact that $N$ is ${\mathfrak A}$-isotropic, we can calculate the normal Jacobi operator $R_N$ from the explicit expression of the curvature tensor of ${Q^*}^{2k}$:
\begin{eqnarray*}
R_N Z & = & R(Z,N)N  \\
& = &  -Z + g(Z,N)N - 3g(Z,\xi)\xi  + g(Z,AN)AN + g(Z,A\xi)A\xi.
\end{eqnarray*}
It follows that $R_N$ has the three constant eigenvalues $0,1,4$ with corresponding eigenbundles $\nu M \oplus ({\mathcal C} \ominus {\mathcal Q}) = \nu M \oplus T_\kappa$, ${\mathcal Q} = T_\lambda \oplus T_\mu$ and ${\mathcal F} = T_\alpha$.
This leads to the following expressions for the Jacobi vector fields along ${\gamma}_{[z]}$:
\begin{equation*}
Z_X(r)=\begin{cases}
(\cosh(2r) - \frac{\alpha}{2}\sinh(2r))E_X(r) & \text{if $X \in T_\alpha$, }\\
(\cosh(r) - \rho\sinh(r))E_X(r) & \text{if $X\in T_\rho$, $\rho\in\{\lambda,\mu\}$,}\\
(1 - \kappa r)E_X(r) & \text{if $X\in T_\kappa$, ${\kappa}{\in}\{0,{\alpha}\}$},
\end{cases}
\end{equation*}
where  ${\alpha}=2\coth(2r)$, ${\lambda}=\coth(r)$, and ${\mu}=\tanh(r)$, and $E_X$ is the parallel vector field along ${\gamma}_{[z]}$ with $E_X(0)=X$.
This shows that ${\rm Ker}(dF) = T_\alpha \oplus T_\lambda$ and thus $\text{rank}F=\text{dim}M-\text{dim}({\rm Ker}dF)=2k$,
where we take into account that in the case $\kappa = \alpha$ the function $1 - 2r\coth(2r)$ is non-vanishing, because if we assume $1-(\coth(r)+\tanh(r))r=0$, then it becomes $r(1+\coth^2(r))=\coth(r)$. From this, taking limit $r{\to}0$, then $0=\lim_{r\to 0}r(1+\coth^2r)=\lim_{r\to 0}\coth(r)=1$. This gives a contradiction for $0<r<1$ and also for $r{\to}1$, and $r(1+{\coth}^2(r)){\not =}\coth(r)$ for $r{\ge}1$.
So, locally, $F$ is a submersion onto a submanifold $H$ of ${Q^*}^{2k}$ of real dimension $2k$. Moreover, the tangent space $T_{F([z])}H$ of $H$ at $F([z])$ is obtained by parallel translation of $(T_\mu \oplus T_\kappa)([z])$ along $\gamma_{[z]}$. Since $T_\mu$ and $T_\kappa$ are both $J$-invariant and $J$ is invariant under parallel translation along geodesics, because ${\nabla}_XJ=0$. Then it follows that the submanifold $H$ is a complex submanifold of ${Q^*}^{2k}$ of complex dimension $k$.

The vector $\eta_{[z]} = \dot{\gamma}_{[z]}(r)$ is a unit normal vector of $H$
at $F([z])$ and the shape operator $S_{\eta_{[z]}}$ of $H$ with respect to
$\eta_{[z]}$ can be calculated from the equation
\[
S_{\eta_{[z]}}Z_X(r) = - Z_X^\prime(r),
\]
where $X \in (T_\mu \oplus T_\kappa)([z])$. The above expression for the Jacobi vector fields $Z_X$ implies $Z_X^\prime(r) = 0$ for $X \in T_\mu([z])$ and $X \in T_{\kappa=0}([z])$, and therefore $S_{\eta_{[z]}} = 0$ if $\kappa = 0$. If $\kappa = \alpha$ we have $S_{\eta_{[z]}}E_X(r) = \frac{2\coth(2r)}{1-2r\coth(2r)} E_X(r)$. Since every complex submanifold of a K\"{a}hler manifold is minimal we must have $2\coth(2r) = 0$. This implies that $(e^{2r})^2+1=0$, which gives a contradiction. So this case ${\kappa}={\alpha}$ can not be occured. Summing up above two cases, we can only consider for ${\kappa}=0$. Then by above theorem, $M$ is locally congruent to an open part of a horosphere or $g(BX,Y)=0$ for even complex $2k$-dimensional complex quadric ${Q^*}^{2k}$.

Now let us consider the next case.
The vectors of the form $\eta_{[q]}$, $[q] \in F^{-1}(\{[z]\})$,
form an open subset of the unit
sphere in the normal space of $H$ at $F([z])$. Since $B_{\eta_{[q]}}$ vanishes
for all $\eta_{[q]}$ it follows that $H$ is a $k$-dimensional totally geodesic complex submanifold of ${Q^*}^{2k}$.
Rigidity of totally geodesic submanifolds now
implies that the entire submanifold
$M$ is an open part of a tube around a $k$-dimensional connected, complete,
totally geodesic complex submanifold $H$ of ${Q^*}^{2k}$.

By the classification of the totally geodesic submanifolds in complex quadrics due to Klein \cite{K}, we can give a corresponding classification for a complex hyperbolic quadric. According to the Klein's classification the submanifold $H$ is either a totally geodesic ${Q^*}^k \subset {Q^*}^{2k}$ or a totally geodesic ${\mathbb C}H^k \subset {Q^*}^{2k}$. The normal spaces of ${Q^*}^k$ are Lie triple systems and the corresponding totally geodesic submanifolds of ${Q^*}^{2k}$ are again $k$-dimensional quadrics. Since $k \geq 2$ it follows that the normal spaces of ${Q^*}^k$ contain all types of tangent vectors of ${Q^*}^{2k}$. This implies that the normal bundle of the tubes around ${Q^*}^k$ contains regular and singular tangent vectors of ${Q^*}^{2k}$. Since the normal bundle of $M$ consists of ${\mathfrak A}$-isotropic tangent vectors only we conclude that $H$ is congruent to ${\mathbb C}H^k$. It follows that $M$ is congruent to an open part of a tube around ${\mathbb C}H^k$. This concludes the proof of Theorem \ref{mainresult}.

\begin{rem}
In Hermitian symmetric space of rank $2$ in compact type, many geometers (see \cite{K}, \cite{R}, \cite{S}, \cite{PS1}, \cite{PS2} and \cite{PSW}) have investigated some geometric structures in compact two-plane Grassmannians, and in non-compact case also many geometers have studied geometric properties in complex hyperbolic two-plane Grassmannians (see \cite{E}, \cite{MR}, \cite{S5} and \cite{S6}). As real Grassmannian manifolds of rank $2$ among the class of Hermitian symmetric spaces we can give two examples of compact complex quadric $Q^m$ and its non-compact dual ${Q^*}^m$. Naturally in complex hyperbolic quadric ${Q^*}^m$ we can consider the
notions of Reeb invariant, Reeb parallel shape operator and parallel Ricci tensor in \cite{S7}, \cite{S8} and \cite{S9} respectively.
\end{rem}

\end{document}